\documentclass[xypic,amscd,syntonly,amssymb,verbatim,12pt]{amsart}

\usepackage{graphicx}
\usepackage[all]{xy}
\usepackage{amssymb}
\usepackage{amsmath}
\usepackage[mathscr]{euscript}

\usepackage{epsfig}
\theoremstyle{plain}
\newtheorem{Thm}{Theorem}
\newtheorem{Prop}[Thm]{Proposition}
\newtheorem{Cor}[Thm]{Corollary}
\newtheorem{Lem}[Thm]{Lemma}

 \theoremstyle{definition}
\theoremstyle{remark}

\numberwithin{equation}{section}

\begin{document}
 %\title{Coupling class of actions of reductive groups}
 %\title{Vertex operators between general $B$-branes}
  \title{Spanning Class in the Category of Branes}

 \author{ ANDR\'{E}S   VI\~{N}A}
\address{Departamento de F\'{i}sica. Universidad de Oviedo.   Garc\'{\i}a Lorca 18.
     33007 Oviedo. Spain. }
 \email{andresvinaescalar@gmail.com}
%\thanks{This work has been partially supported by Ministerio de Ciencia y
%Tecnolog\'{\i}a, grant FPA2009-11061}

 \maketitle
\begin{abstract}

Given a generic anticanonical hypersurface $Y$ of a toric variety determined by a reflexive polytope, we define a line bundle ${\mathcal L}$  on $Y$ that generates a
%such that ${\{\mathcal L}^{\otimes i}\}$ is
 spanning class in the bounded derivative category $D^b(Y)$. From this fact,
 %fact that ${\mathcal L}$ is a spanning class, 
we deduce properties of some spaces of strings    related with the brane ${\mathcal L}$. We prove a vanishing theorem for the 
%spaces of
 vertex operators associated to strings stretching  from  branes of the form ${\mathcal L}^{\otimes i}$ to nonzero objects in $D^b(Y)$.  We also define a gauge field on ${\mathcal L}$ which minimizes  the corresponding Yang-Mills functional.

\end{abstract}

\section {Introduction} \label{S:intro}

%A well-known fact is 
 It is well-known that a compact toric manifold is not Calabi-Yau.
 % However, certain hypersurfaces in some toric varieties are Calabi-Yau. More precisely,
 However,  Batyrev \cite{Batyrev} showed the existence of anticanonical hypersurfaces,  in    the toric variety $X$ determined by a reflexive polytope,   that are Calabi-Yau. 
In this note, we will prove some particular properties of $D$-branes, strings,  vertex operators and gauge fields on these hypersurfaces. 

 We will denote by $Y$ such a hypersurface.   The $D$-branes of type $B$ on $Y$ are the objects of $D^b(Y)$, the   bounded derived category of coherent sheaves on $Y$ (see monograph \cite{Aspin-et}). In this context, 
 given the branes ${\mathcal F}_1$ and ${\mathcal F}_2$, the strings with ghost number $k$ stretching from the brane ${\mathcal F}_1$ to the brane ${\mathcal F}_2$ are the elements of the ext group ${\mathit Ext}^k({\mathcal F}_1,\,{\mathcal F}_2)$  \cite[Section 5.2]{Aspin}.
 
Let $\Delta$ be a reflexive polytope in $({\mathbb R}^n)^*$ that defines the toric variety $X$.  The multiple $(n-1)(-K_X)$ of the anticanonical divisor of $X$ determines a line bundle on $X$, and we put ${\mathcal L}$  for the restriction of that bundle to $Y$.
In Propositions \ref{Prop:Spaning} and \ref{Prop:Spaning1},  we will prove 
 %the existence of a line bundle ${\mathcal L}$ on $Y$ such 
that the family $\{{\mathcal L}^{\otimes i}\,|\,i\in{\mathbb Z}\}$ is a spanning class of $D^b(Y)$ \cite{Bridgeland, Huybrechts}.
 Thus, given ${\mathcal F}$   a nonzero  brane on $Y$,   there  exist two integers $a,r$ such that the space of strings 
\begin{equation}\label{Eq:Exta}
{\mathit Ext}^r({\mathcal L}^{\otimes a},\,{\mathcal F})\ne 0.
\end{equation}

Condition (\ref{Eq:Exta}) implies relations between certain spaces of strings. If ${\mathcal H}$ is a brane such that can bind with the above brane ${\mathcal L}^{\otimes a}$ to form a new brane ${\mathcal G}$, then we deduce relations between the spaces of strings
 ${\mathit Ext}^j({\mathcal F},\,{\mathcal G})$ and  ${\mathit Ext}^j({\mathcal F},\,{\mathcal H})$,  which are gathered in Theorem \ref{Thm:Ext1}.

In the case that a brane ${\mathcal G}$ can decay to the above brane ${\mathcal F}$ and  other one ${\mathcal J}$, then we prove relations between the spaces 
  ${\mathit Ext}^j({\mathcal J},\,{\mathcal L}^{\otimes a})$ and  ${\mathit Ext}^j({\mathcal G},\,{\mathcal L}^{\otimes a})$, which we have collected in Theorem \ref{Thm:Ext2}.

	 We set ${\mathcal Hom}(\,.\,,\,.\,) $ for the sheaf functor Hom (see \cite[page 87]{{Kas-Sch}}) of the category $\mathfrak{Coh}(Y)$ of coherent sheaves on $Y$
$$ {\mathcal Hom}(\,.\,,\,.\,):\mathfrak{Coh}(Y)^{\rm op}\times\mathfrak{Coh}(Y)\to\mathfrak{Sh},$$
 where $\mathfrak{Sh}$ is the category of  sheaves of ${\mathbb C}$-vector spaces  on $Y$.
% It is easy to check that
 % \begin{equation}\label{OHom}
 % {\mathcal Hom}\big({\mathcal O}(V_1),\,{\mathcal O} (V_2)  \big)={\mathcal O}({\rm Hom}(V_1,\,V_2)   ).
 %  \end{equation}
	%As we said, the bounded derived category of $\mathfrak{Coh}(X)$ will be denoted by $D(X)$, 
	Thus, one has the derived functor
 $$R{\mathcal Hom}(\,.\,,\,.\,): D^b(Y)^{\rm op}\times D^b(Y)\to D(\mathfrak{Sh}),$$
  where $D(\mathfrak{Sh})$ is the derived category of $\mathfrak{Sh}$.
	%and $D(Y)$ the derived category $\mathfrak{Coh}(Y)$. 
	By definition
  ${\mathcal Ext}^k({\mathcal F},\,{\mathcal G})=  H^k R{\mathcal Hom}({\mathcal F},\,{\mathcal G}).$
	
	The space of vertex operators for strings, with ghost number $k$, between the branes ${\mathcal F}$ and ${\mathcal G}$ is \cite{Vina}
	$$\bigoplus_qH^q(Y,\,{\mathcal Ext}^k({\mathcal F},\,{\mathcal G})).$$ 
	
	We determine a  rectangle ${\sf R}$ in ${\mathbb R}^2$ where ``are located''
	%the pairs $(q,k)$ corresponding to 
	the nontrivial spaces of vertex operators for strings from ${\mathcal L }^{\otimes b}$ and ${\mathcal F}$, $b$ being an arbitrary integer. That is, if the space of vertex operators $H^q(Y,\,{\mathcal Ext}^p({\mathcal L}^{\otimes b},\,{\mathcal F}))$ is nonzero, then $(p,\,q)$ belongs to ${\sf R}$ (Theorem \ref{Thm:parallelogram}).  From this theorem, it turns out the existence of an
	%There exists an 
	 integer $i_0$, such that if $i\leq i_0$  the   space of vertex operators $H^q(Y,\,{\mathcal Ext}^p({\mathcal L}^{\otimes i},\,{\mathcal F}))$
	corresponding to a some point of ${\sf R}$ 
	is in fact nontrivial (Corollary \ref{Cor:Vertex}).
	
%Let $\Delta$ be the reflexive lattice polytope in $({\mathbb R}^n)^*$ which determine the variety $X$. Then the   $(n-1)\Delta$ is very ample polytope, and its lattice points   define an embedding in a projective space
	% $Y\overset{j}{\hookrightarrow}{\mathbb P}^N$, such that ${\mathcal L} =j^*({\mathcal O}_{{\mathbb P}^N}(1))$.
	Although $Y$ is not smooth, we define in Section \ref{S: Yang-Mills fields} 
	 connections   on ${\mathcal L}$ by using K\"ahler differentials. As $(n-1)\Delta$ is very ample polytope,  its lattice points determine an embedding in a projective space
	 $Y\overset{j}{\hookrightarrow}{\mathbb P}^N$, such that ${\mathcal L} =j^*({\mathcal O}_{{\mathbb P}^N}(1))$.
	%From a connection on ${\mathcal L}$, which is 
	By pulling back  the Chern connection on ${\mathcal O}_{{\mathbb P}^N}(1),$ we define a particular connection  on ${\mathcal L}$.
	  %We determine the expression in coordinates of this gauge field and the one of its curvature.
	   In Proposition \ref{Prop:YM}, we prove that this gauge field minimizes the Yang-Mills action.

In the case that $Y$ is a  hypersurface   invariant under the action of the torus $T$ of the variety $X$, then ${\mathcal L}$ is a $T$-equivariant brane on X \cite{V4}. In Section \ref{S:Equiv}, we determine the value of the localization of  
%$c_1^T({\mathcal L})$, 
the  $T$-equivariant Chern class of ${\mathcal L}$. This result is stated in Theorem \ref{Th:ewquiv}.  

%%%%%%%%%%%%%%%%%%%%%%%%%%%%%%%%%%%%%%%%%%%%%%%%%%%%%%%%%%%%%%%%%%%%%%%%%%%%%%%%%%%%%%%%%%%%%%%%%%%%%%%%%%%%%%%%%%%%%%%%%%%%%%%
%%%%%%%%%%%%%%%%%%%%%%%%%%%%%%%%%%%%%%%%%%%%%%%%%%%%%%%%%%%%%%%%%%%%%%%%%%%%%%%%%%%%%%%%%%%%%%%%%%%%%%%%%%%%%%%%%%%%%%%%%%%%%%%%%%%%%%%%%%%%%%%%%

\section{Calabi-Yau hypersurfaces}\label{S:Calabi}

\noindent
{\bf Reflexive polytopes.} Let $\Delta$ be a lattice polytope of dimension $n$ in $({\mathbb R}^n)^*$. Given a facet $F$ of $\Delta$, there is   a unique vector $v_F\in{\mathbb Z}^n$ conornal to $F$ and inward to $\Delta$. So, $F$ is on  the hyperplane in $({\mathbb R}^n)^*$ of equation
 $\langle m,\,v_F\rangle=-c_F$, with $c_F\in {\mathbb Z}$ and
$$\Delta=\bigcap_{F\, {\rm facet}}\{  m\in ({\mathbb R}^n)^*\,|\, \langle m,\,v_F\rangle\geq -c_F   \}.$$

We denote by $X$ the toric variety determined by $\Delta$ and by $D_F$ the divisor on $X$ associated to the facet $F$. The divisor $D_{\Delta}:=\sum_F c_FD_F$ is Cartier, ample and base pointfree \cite[page 269]{C-L-S}. Moreover, the torus invariant divisor $K_X:=-\sum_F D_F$ is a canonical divisor of $X$.   

From now on, we assume that $\Delta$ is a {\it reflexive} polytope; that is,   $c_F=1$ for all facet $F$. Then the divisor $K_X$ is Cartier and the variety $X$ is Gorenstein. Thus, the canonical sheaf $\omega_X={\mathcal O}_X(K_X)$ is a line bundle and $X$ is Fano (that is, the anticanonical divisor $-K_X$ is ample).

On the other hand, a generic element $Y$ of the linear system $|-K_X|$ is an orbifold. By the adjunction formula,  the canonical sheaf of $Y$, $\omega_Y$,  satisfies
$$\omega_Y\simeq\omega_X\otimes_{{\mathcal O}_X}({\mathcal I}_Y/{\mathcal I}_Y^2)^{\vee},$$ 
where ${\mathcal I}_Y$ denotes the ideal sheaf of $Y$. As $\omega_X={\mathcal O}_X(-Y)$ and ${\mathcal I}_Y/{\mathcal I}_Y^2={\mathcal O}_X(-Y)\otimes_{{\mathcal O}_X} {\mathcal O}_Y,$ ones deduces the triviality of $\omega_Y$. The hypersurface $Y$ is in fact a Calabi-Yau variety (for more details see \cite[Sec. 4.1]{C-K}).

\smallskip

\noindent
{\bf Homological dimension.} Henceforth, $Y$ will be a generic anticanonical hypersurface of $X$.
%By $D^b(Y)$ we denote the bounded derived category of coherent sheaves on $Y$.
 Since $Y$ is a Gorenstein variety a Serre functor 
%${\sf S}$  
$${\sf S}:D^b(Y)\to D^b(Y),$$
is ${\sf S}=(\,.\,)\otimes \omega_Y[n-1]$, where $[n-1]$ denotes the shifting of the complex by $n-1$ to the left \cite{B-K}. By the Serre duality, we have functorial isomorphisms
$${\rm Hom}_{D^b(Y)}({\mathcal F},\,{\mathcal G})\simeq {\rm Hom}_{D^b(Y)}({\mathcal G},\,{\sf S}({\mathcal F}))^{\vee}$$
for every objects ${\mathcal F}$, ${\mathcal G}$ of $D^b(Y)$.

In particular,   if ${\mathcal A}$, ${\mathcal B}$ are objects of the category $\mathfrak{Coh}(Y)$ of coherent sheaves on $Y$, then
$${\mathit Ext}^k({\mathcal A},\,{\mathcal B})= {\rm Hom}_{D^b(Y)}({\mathcal A},\,{\mathcal B}[k])\simeq\mathit{  Ext}^{n-1-k}({\mathcal B},\,{\mathcal A}\otimes\omega_Y)^{\vee}.$$
 As ${\mathit Ext}^{n-1-k}({\mathcal B},\,{\mathcal A}\otimes\omega_Y)=0$, for $n-1<k$, the homological dimension of the category $\mathfrak{Coh}(Y)$  is $\leq n-1$. On the other hand, taking ${\mathcal A}={\mathcal O}_Y$ and ${\mathcal B}=\omega_Y$, one has
 $${\mathit Ext}^{n-1}({\mathcal O}_Y,\,\omega_Y  )=\mathit{Hom}(\omega_Y,\,\omega_Y)^{\vee}=\mathit{Hom}({\mathcal O}_Y,\,{\mathcal O}_Y)^{\vee}.$$
 As this space is different from zero,
 we deduce that the mentioned homological dimension 
 %of category $\mathfrak{Coh}(Y)$ 
 is $n-1$. 
 Therefore, given ${\mathcal B}$ an object of $\mathfrak{Coh}(Y)$, there exists $p_0$, such that 
 $${\rm Hom}_{D^b(Y)}({\mathcal A} ,\,{\mathcal B}[p])=0,$$
  for all $p>p_0$ and for any object ${\mathcal A}$ of $\mathfrak{Coh}(Y)$. 
	
	Let ${\mathcal F}^{\bullet}$ be an object in the  bounded derived category $D^b(Y)$, then 
 $$q_0:={\rm max}\{m\,|\, H^{-m}({\mathcal F}^{\bullet})\ne 0   \}<\infty.$$
 On the other hand, the Grothendieck spectral sequence
 $$E_2^{pq}:={\rm Hom}_{D^b(Y)}\big( H^{-q}({\mathcal F}^{\bullet}),\,{\mathcal B}[p]    \big)$$
 abuts to ${\rm Hom}_{D^b(Y)}({\mathcal F}^{\bullet},\,{\mathcal B}[p+q])$. As $E_2^{pq}=0$ whether $p>p_0$ or $q>q_0$, we deduce the following lemma.
 
 \begin{Lem} \label{Lem:spectral}
 Given an object ${\mathcal F}^{\bullet}$ of $D^b(Y)$ and a coherent sheaf ${\mathcal B}$ on $Y$, there exists an integer $j_0$ such that ${\rm Hom}_{D^b(Y)}({\mathcal F}^{\bullet},\,{\mathcal B}[j])=0$ for all $j>j_0$.
 \end{Lem}
 An analog reasoning proves the existence of an integer $j_1$ such that ${\rm Hom}_{D^b(Y)}({\mathcal F}^{\bullet},\,{\mathcal B}[j])=0$ for all $j<j_1$.

 \smallskip

\noindent
{\bf Spanning class.}
The multiple $-(n-1)K_X$ of the anticanonical divisor $-K_X$ is a very ample divisor on $X$ \cite[page 71]{C-L-S}. We put ${\mathcal L}$ for the restriction to $Y$ of the very ample invertible sheaf on $X$ defined by the Cartier divisor $-(n-1)K_X$.

Let ${\mathcal F}^{\bullet}$ be a nonzero object of the bounded derived category $D^b(Y)$. Then there exists an integer $r$, such that the cohomology $H^r({\mathcal F}^{\bullet})\ne 0$ and $H^j({\mathcal F}^{\bullet})= 0$, for $j<r$. Thus, we may assume that ${\mathcal F}^{\bullet}$ is a complex of the form
$$\dots\to 0\to 0\to{\mathcal F}^r\overset{\partial^r}{\rightarrow} {\mathcal F}^{r+1}\to\dots.$$
As $H^r({\mathcal F}^{\bullet})={\rm Ker}(\partial^r)\subset{\mathcal F}^r$, we have a monomorphism from the complex
 $$\dots\to 0\to 0\to H^r({\mathcal F}^{\bullet})\to 0\to 0\dots$$ 
(with only one nonzero element at the position $0$) to the complex ${\mathcal F}^{\bullet}[r]$.

Since the natural functor $\mathfrak{Coh}(Y)\to D(Y)$ is fully faithful \cite[page 165]{Ge-Ma},  we have
the following injective map, induced by the above monomorphism,
\begin{align}\label{aligneament} \notag{\rm Hom}_{\mathfrak{Coh}(Y)}({\mathcal L}^{\otimes i},\, H^r({\mathcal F}^{\bullet})) =
&{\rm Hom}_{D^b(Y)}({\mathcal L}^{\otimes i},\, H^r({\mathcal F}^{\bullet}))\to \\ 
 &{\rm Hom}_{D^b(Y)}({\mathcal L}^{\otimes i},\, {\mathcal F}^{\bullet}[r]).
\end{align}

As ${\mathcal L}$ is very ample and $H^r({\mathcal F}^{\bullet})$ is a coherent sheaf on $Y$, there exists an integer $m_0$ such that, for all $m\geq m_0$, the sheaf $H^r({\mathcal F}^{\bullet})\otimes_{{\mathcal O}_Y}{\mathcal L}^{\otimes m} $ is generated by a finite number of global sections  \cite[page 121]{Hart}.  The space of global sections
\begin{align}\label{alignea}
&\Gamma \big(Y,\,H^r({\mathcal F}^{\bullet})\otimes_{{\mathcal O}_Y}{\mathcal L}^{\otimes m}\big)=
\Gamma\big(Y,\,{\mathcal Hom}({\mathcal L}^{\otimes -m},\,H^r({\mathcal F}^{\bullet}))\big)=\\ \notag
&{\rm Hom}_{\mathfrak{Coh}(Y)}({\mathcal L}^{\otimes -m},\,H^r({\mathcal F}^{\bullet})).
\end{align}

 If ${\rm Hom}_{D^b(Y)}({\mathcal L}^{\otimes i},\, {\mathcal F}^{\bullet}[r])$ {\em were zero for} for an integer $i\leq -m_0$, then from (\ref{aligneament}) it would follow 
 ${\rm Hom}_{\mathfrak{Coh}(Y)}({\mathcal L}^{\otimes i},\, H^r({\mathcal F}^{\bullet}))=0$, and from (\ref{alignea}) we would deduce that the $ \Gamma \big(Y,\,H^r({\mathcal F}^{\bullet})\otimes_{{\mathcal O}_Y}{\mathcal L}^{\otimes m}\big)=0$, for  $m:=-i\geq m_0$. But this is in contradiction with the fact that  the non trivial sheaf $H^r({\mathcal F}^{\bullet})\otimes_{{\mathcal O}_Y}{\mathcal L}^{\otimes m}$ is generated by its global sections, if $m\geq m_0$. We have proved the following proposition.

\begin{Prop}\label{Prop:Spaning} If ${\mathcal F}$ is a nonzero object of the bounded derived category $D^b(Y)$, then there 
exist integers $i_0:=i_0({\mathcal F})$ and $r=r({\mathcal F})$ such that
${\mathit Ext}^r({\mathcal L}^{\otimes i},\,{\mathcal F})\ne 0$, for any $i\leq i_0$.
\end{Prop}

\begin{Prop}\label{Prop:Spaning1} If ${\mathcal G}$ is a nonzero object of the bounded derived category $D^b(Y)$, then there exist integers $n_0,l$ such that
${\mathit Ext}^l({\mathcal G},\;  {\mathcal L}^{\otimes i})\ne 0$, for any $i\leq n_0$.
\end{Prop}

{\it Proof.} By  the Serre duality, for all $l,i$
\begin{align} \notag
\mathit{Ext}^l({\mathcal G},\;  {\mathcal L}^{\otimes i})=&{\rm Hom}_{D^b(Y)}\big({\mathcal G},\, {\mathcal L}^{\otimes i}[l]\big)\simeq{\rm Hom}_{D^b(Y)} \big({\mathcal L}^{\otimes i}[l],\,{\sf S}({\mathcal G} ) \big)^{\vee} \\ \notag
\simeq & {\rm Hom}_{D^b(Y)}\big({\mathcal L}^{\otimes i},\,{\sf S}({\mathcal G} )[-l] \big)^{\vee}\simeq
 \mathit{Ext}^{-l}({\mathcal L}^{\otimes i},\,{\sf S}( {\mathcal G}))^{\vee}.
\end{align}
As the Serre functor ${\sf S}$ is an equivalence, ${\sf S}( {\mathcal G})\ne 0$ and applying to it Proposition \ref{Prop:Spaning}, we can  take $n_0= i_0({\sf S}( {\mathcal G}))$ and $l=-r({\sf S}( {\mathcal G}))$.

\qed

%\begin{Cor}\label{Cor:Spaning}
%If ${\mathcal G}$ is an object of $D^b(Y)$ such that ${\mathit Ext}^r(\mathcal{G},\,{\mathcal L}^{\otimes i})= 0$, for all integers $i,r$, then ${\mathcal G}=0$. 
%\end{Cor}
%{\it Proof.} By the hypothesis and the Serre duality, for all $r,i$
%\begin{align} \notag
%0=&{\rm Hom}_{D^b(Y)}\big({\mathcal G},\, {\mathcal L}^{\otimes i}[r]\big)\simeq{\rm Hom}_{D^b(Y)} \big({\mathcal L}^{\otimes i}[r],\,{\sf S}({\mathcal G} ) \big)^* \\ \notag
%\simeq & {\rm Hom}_{D^b(Y)}\big({\mathcal L}^{\otimes i},\,{\sf S}({\mathcal G} )[r] \big)^*.
%\end{align}
%By the proposition, ${\sf S}({\mathcal G} )=0$. As  the Serre functor ${\sf S}$ is an equivalence, we conclude that ${\mathcal G}=0$.
%\qed
 
In summary, the family of powers $\{{\mathcal L}^{\otimes a}\,|\,a\in{\mathbb Z} \}$ is a spanning class \cite{Bridgeland}  in $D^b(Y)$. That is, the vanishing 
 ${\mathit Ext}^k({\mathcal F},\, {\mathcal L}^{\otimes a})=0$ for all $a$ and all $k$ implies ${\mathcal F}=0$. Equivalently, if
${\mathit Ext}^k( {\mathcal L}^{\otimes a},\,{\mathcal F})=0$ for all $a$ and all $k$, then ${\mathcal F}=0$.

 %%%%%%%%%%%%%%%%%%%%%%%%%%%%%%%%%%%%%%%%%%%%%%%%%%%%%%%%%%%%%%%%%%%%%%%%%%%%%%%%%%%%%%%%%%%%%%%%%%%%%%%%%%%%%%%%%%%%%%%%%%%%%%%%%%%%%%%%%%%%
%%%%%%%%%%%%%%%%%%%%%%%%%%%%%%%%%%%%%%%%%%%%%%%%%%%%%%%%%%%%%%%%%%%%%%%%%%%%%%%%%%%%%%%%%%%%%%%%%%%%%%%%%%%%%%%%%%%%%%%%%%%%%%%%%%%%%%%%%%%%%%%%%%%%%%%

\section {$B$-branes on the hypersurface} 

As it was mentioned in the Introduction, the branes of type $B$ on the Calabi-Yau orbifold $Y$ are the objects of  bounded derived category $D^b(Y)$, and   the space of strings with ghost number $k$ stretching between the branes ${\mathcal F}_1$ and ${\mathcal F}_2$ is the group ${\mathit Ext}^k({\mathcal F}_1,\, {\mathcal F}_2)$.

 Given a nonzero brane ${\mathcal F}$, by Proposition \ref{Prop:Spaning1},  there exists $a\in{\mathbb Z}$ such that the following set is nonempty 
 $$S:=\{r\in{\mathbb Z},\,|\,\mathit{  Ext}^r({\mathcal F},\,{\mathcal L}^{\otimes a})\ne 0\}\ne\emptyset.$$
 By Lemma \ref{Lem:spectral} together with the comment after that lemma, $S$ is finite, and we denote by $k_1$ and $k_2$ the minimum and the maximum of $S$, respectively.
 
 Let us assume that  the branes ${\mathcal L}^{\otimes a}$ and ${\mathcal H}$ may bind to form the brane ${\mathcal G}$ through a potentially tachyonic open string. In mathematical terms, one has a distinguished triangle 
 \begin{equation}\label{DistTrian}
 {\mathcal L}^{\otimes a} \to{\mathcal G}\to {\mathcal H}\overset{+1}{\longrightarrow}
 \end{equation}
 in the category $D^b(Y)$. We have the corresponding long exact Ext sequence
 $$\to{\mathit Ext}^j( {\mathcal F},\,{\mathcal L}^{\otimes a}  )\to{\mathit Ext}^j( {\mathcal F},\,{\mathcal G}  )\to
{\mathit Ext}^j( {\mathcal F},\,{\mathcal H} )\to
{\mathit Ext}^{j+1}( {\mathcal F},\,{\mathcal L}^{\otimes a}  )\to 
$$
By definition of $k_1$ and $k_2$, one deduces the following theorem.
\begin{Thm}\label{Thm:Ext1} With the notations introduced above,
\begin{enumerate}
\item Whether $j<k_1-1$ or $k_2<j$
$$\mathit{  Ext}^j( {\mathcal F},\,{\mathcal G}  )\simeq
{\mathit Ext}^j( {\mathcal F},\,{\mathcal H} ).$$
\item   For $j\in\{k_1,\,k_2-1\},$ the spaces ${\mathit  Ext}^j( {\mathcal F},\,{\mathcal G}  )$ and
${\mathit Ext}^j( {\mathcal F},\,{\mathcal H} )$ are not isomorphic.
 \item For $j=k_1-1$, each string with ghost number $j$ from ${\mathcal F}$ to ${\mathcal G}$ admits a unique extension to a string 
between ${\mathcal F}$ and ${\mathcal H}$.
\item For $j=k_2$, each string with ghost number $j$ stretching from ${\mathcal F} $ to ${\mathcal H}$ 
%can be lifted to a string 
is the extension of a string from ${\mathcal F}$ to ${\mathcal G}$. 
\end{enumerate}
\end{Thm}
 In other words, for the numbers $j$ mentioned in the first item  of theorem, each string with ghost number $j$ from ${\mathcal F}$ to ${\mathcal H}$ can be uniquely lifted  to a string  from ${\mathcal F}$ to ${\mathcal G}$. 

\smallskip

On the other hand, if there is a string between the branes ${\mathcal F}$ and ${\mathcal G}$, then ${\mathcal G}$ can decay to ${\mathcal F}$ and other brane ${\mathcal J}$. We have the corresponding distinguished triangle
$${\mathcal F}\to{\mathcal G}\to{\mathcal J}\overset{+1}{\longrightarrow},$$
and the respective long exact sequence of Ext groups
$$\to{\mathit Ext}^{j-1}( {\mathcal F},\,{\mathcal L}^{\otimes a}  )\to\mathit{  Ext}^{j}( {\mathcal J},\, {\mathcal L}^{\otimes a} )\to
{\mathit Ext}^{j}( {\mathcal G},\,{\mathcal L}^{\otimes a} )\to
{\mathit Ext}^{j}( {\mathcal F},\,{\mathcal L}^{\otimes a}  )\to
$$
\begin{Thm}\label{Thm:Ext2} With the above notations, 
 \begin{enumerate} 
\item Whether $j<k_1$ or $k_2+1<j$
$$\mathit{  Ext}^j( {\mathcal J},\, {\mathcal L}^{\otimes a}  )\simeq
{\mathit Ext}^j( {\mathcal G},\, {\mathcal L}^{\otimes a} ).$$
\item For $j\in\{k_1+1,\,k_2\}$, $\mathit{  Ext}^j( {\mathcal J},\, {\mathcal L}^{\otimes a}  )$ and
${\mathit Ext}^j( {\mathcal G},\, {\mathcal L}^{\otimes a})$ are not isomorphic.  
\item For $j=k_1$, each string with ghost number $j$ between ${\mathcal J} $ to ${\mathcal L}^{\otimes a}$ admits a unique  lift of a string from ${\mathcal G}$ to ${\mathcal L}^{\otimes a}$. 
\item For $j=k_2+1$,  each string with ghost number $j$ from ${\mathcal G}$ to ${\mathcal L}^{\otimes a}$ is the lift of a string   
between ${\mathcal J}$ and ${\mathcal L}^{\otimes a}$.
\end{enumerate}
\end{Thm}
That is, for the $j$'s considered in the first item of theorem, each string with ghost number $j$ stretching
 from ${\mathcal G}$ to ${\mathcal L}^{\otimes a}$ can be  uniquely  extended   to a string from
  ${\mathcal J}$ to ${\mathcal L}^{\otimes a}$.

%%%%%%%%%%%%%%%%%%%%%%%%%%%%%%%%%%%%%%%%%%%%%%%%%%%%%%%%%%%%%%%%%%%%%%%%%%%%%%%%%%%%%%%%%%%%%%%%%%%%%%%%%%%%%%%%%%%%%%%%%%%%%%%%%%%%%%%%%%%%%%%%%%%%%%%%%%%%%%%%%%%
%%%%%%%%%%%%%%%%%%%%%%%%%%%%%%%%%%%%%%%%%%%%%%%%%%%%%%%%%%%%%%%%%%%%%%%%%%%%%%%%%%%%%%%%%%%%%%%%%%%%%%%%%%%%%%%%%%%%%%%%%%%%%%%%%%%%%%%%%%%%%%%%%%%%%%%%%%%%%%%%

\section {Vertex operators}

Let ${\mathcal F}^{\bullet}$ be a nonzero object of the bounded derived category $D^b(Y)$. Then there exist two integers $r,\,s$ such that
\begin{equation}\label{Hrmathcal}
H^r({\mathcal F}^{\bullet})\ne 0\ne H^s({\mathcal F}^{\bullet}),\;\;\hbox{and }\; H^j({\mathcal F}^{\bullet})=0\;\; \hbox{for}\;j\notin[r,\,s].
\end{equation}
In other words, ${\mathcal F}^{\bullet}$ is an object of the category $ D^{\leq s}\cap  D^{\geq r}.$ 

Given      $b\in\mathbb Z$, we put ${\mathcal A}:={\mathcal L}^{\otimes b}$. The spectral sequence
$$E^{p,q}_2={\rm Hom}_{D^b(Y)}({\mathcal A},\,H^q({\mathcal F}^{\bullet})[p]),$$
converges to
$${\rm Hom}_{D^b(Y)}({\mathcal A},\,{\mathcal F}^{\bullet}[p+q])= {\mathit Ext}^{p+q}({\mathcal A},\,{\mathcal F}^{\bullet}) .$$
Since 
%$E_2^{p,q}=0$ 
$H^q({\mathcal F}^{\bullet})[p]=0$
for all $(p,\,q)$ such that $p+q\notin[r,s]$, it follows that 
%$E^{p,q}_{\infty}=0$, when $p+q\notin[r,\,s]$. That is,
${\mathit Ext}^k({\mathcal A},\,{\mathcal F}^{\bullet})=0$, for $k\notin[r,\,s]$.
Thus, the result stated in Proposition \ref{Prop:Spaning} can be improved slightly. 
% admits a slightly stronger version.
\begin{Prop}\label{Prop:strong}
 Let ${\mathcal F}$ be above   nonzero brane on $Y$. Then there exists an integer $i_0$ such that, if $i\leq i_0$, for some ghost number $k$, with $r\leq k \leq s$,    the corresponding space of strings ${\mathit Ext}^k({\mathcal L}^{\otimes i},\,{\mathcal F})$ is nonzero. 
 \end{Prop}

As ${\mathcal F}^{\bullet}$ satisfies  (\ref{Hrmathcal}), there exists an complex ${\mathcal G}^{\bullet}$ quasi-isomorphic to  
${\mathcal F}^{\bullet}$, such that $ {\mathcal G}^{j}=0$ for $j\notin [r,s]$. Thus  ${\mathcal F}^{\bullet}$ and ${\mathcal G}^{\bullet}$ are isomorphic objects in the derived category $D^b(Y)$. In particular 
 $${\mathcal Ext}^p({\mathcal A},\,{\mathcal F}^{\bullet})={\mathcal Ext}^p({\mathcal A},\,{\mathcal G}^{\bullet}).$$

Since ${\mathcal A}$ is a locally free ${\mathcal O}_Y$-module,
$${\mathcal A}_{\bullet}:\;\;\;\to 0\to 0\to {\mathcal A}\overset{1}{\to}{\mathcal A}$$
is a locally free resolution of ${\mathcal A }.$ 
By definition  of the Ext sheaves
$${\mathcal Ext}^p({\mathcal A},\,{\mathcal G}^{\bullet})=H^p({\mathcal Hom}^{\bullet}({\mathcal A}_{\bullet},\,{\mathcal G}^{\bullet})),$$
where
$${\mathcal Hom}^k({\mathcal A}_{\bullet},\,{\mathcal G}^{\bullet})=\prod_i{\mathcal Hom}({\mathcal A}_i,\,{\mathcal G}^{i+k})={\mathcal Hom}({\mathcal A},{\mathcal G}^k).$$
 Hence, 
 \begin{equation}\label{Ext1}
 {\mathcal Ext}^p({\mathcal A},\,{\mathcal F}^{\bullet})=0,\;\;\hbox{for}\, p\notin[r,\,s]
 \end{equation}

By  Grothendieck's vanishing theorem $H^q(Y,\,{\mathcal Ext}^p({\mathcal A},\,{\mathcal F}^{\bullet}))=0$, for $q\notin [0,\,n-1]$. 
We have proved the following theorem.
\begin{Thm}\label{Thm:parallelogram}
If the space of vertex operators $H^q(Y,\,{\mathcal Ext}^p({\mathcal L}^{\otimes b},\,{\mathcal F}))$ for strings between ${\mathcal L}^{\otimes b}$ and ${\mathcal F}\in D^{\leq s}\cap  D^{\geq r}$ is nonzero, then $(p,\,q)$ is a point of the rectangle $[r,\,s]\times[0,\,n-1]$.
% parallelogram ${\sf P}$.
\end{Thm}

Since ${\mathcal L}$ generates a spanning class in $D^b(Y)$, there are integer $a$ and $c$, such that
 ${\mathit Ext}^c({\mathcal L}^{\otimes a},\,{\mathcal F}^{\bullet})\ne 0$.
From the local-to-global spectral sequence
$$H^q(Y,\,{\mathcal Ext}^p({\mathcal L}^{\otimes a},\,{\mathcal F}^{\bullet}))\Longrightarrow {\mathit Ext}^{p+q}({\mathcal L}^{\otimes a},\,{\mathcal F}^{\bullet}),$$
%together with Theorem \ref{Thm:parallelogram},
 we conclude that the following space of vertex operators 
$$\bigoplus_{p+q=c} H^q(Y,\,{\mathcal Ext}^p({\mathcal L}^{\otimes a},\,{\mathcal F}^{\bullet}))\ne 0.$$
%${\mathit Ext}^{c}({\mathcal A},\,{\mathcal F}^{\bullet})=0$, for $c\notin[r,\,n-1+s]$.

Thus, Proposition \ref{Prop:Spaning} gives the following corollary.
\begin{Cor} \label{Cor:Vertex}
 Given a nonzero object ${\mathcal F}$ of $ D^{\leq s}\cap  D^{\geq r}$, there exists an integer $i_0$ satisfying the following property:
For each $i\leq i_0$ there exists a point $(p,\,q)\in[r,\,s]\times[0,\,n-1]$ such that the vertex space $H^q(Y,\,{\mathcal Ext}^p({\mathcal L}^{\otimes i},\,{\mathcal F}))$ is nontrivial.
\end{Cor}

%%%%%%%%%%%%%%%%%%%%%%%%%%%%%%%%%%%%%%%%%%%%%%%%%%%%%%%%%%%%%%%%%%%%%%%%%%%%%%%%%%%%%%%%%%%%%%%%%%%%%%%%%%%%%%%%%%%%%%%%%%%%%%%%%%%%%%%%%%%%%%%%%%%%%%%%%%%%%%%%%%%%%%%%%%%%%%%%%%%%%%%%%%%%%%%%%%%%%%%%%%%%%%%%%%%%%%%%%%%%%%%%%%%%%%%%%%%%%%%%%%%%%%%%%%%%%%%%%%%%%%%%%%%%%%%%%%%%%%%%%%%%%%%%

\section{Yang-Mills fields on $Y$
% ${\mathcal L}$
 }\label{S: Yang-Mills fields}
 
 \noindent
 {\bf Embedding of $Y$ in ${\mathbb P}^N$.}
%Obviously, the toric varieties defined by the polytopes $\Delta$ and $(n-1)\Delta$ are equal.
We denote by $N$ the number of lattice points of $(n-1)\Delta$.
%; i. e.  $ N:=\#\big((n-1)\Delta\cap {\mathbb Z}^n\big)$.
 Each $m=(a_1,\dots, a_n)\in ((n-1)\Delta)\cap {\mathbb Z}^n$ determines a character of the torus $T$ of the variety $X,$ 
$$\chi_{m}: T=({\mathbb C}^{\times})^n\to{\mathbb C}^{\times},\;\;\chi_{m}(z_1,\dots,z_n)=\prod_{i=1}^n z_i^{a_i}.$$

 Since the divisor $D=(n-1)\sum_FD_F$ is very ample \cite[page 269]{C-L-S},
%The divisor $D_{\Delta}$, defined in Section \ref{S:Calabi},  $D_{\Delta}=\sum_FD_F=-K_X$ is Cartier and ample; moreover,  the divisor
 %$D:=(n-1)D_{\Delta}$ i{}s very ample \cite[page 269]{C-L-S}.
%  But the divisor $D(-{}(n-1)K_X)$ associated to $(n-1)(-K_X)$ is very ample. 
we have the embedding $\varphi:X\to{\mathbb P}^N$ determined by the map
%$$\varphi:X\to{\mathbb P}^N,\;\; \varphi(x)=(\chi_{m_1}(x):\dots :\chi_{m_N}(x)),$$
$$x\in T\mapsto (\chi_{m_1}(x):\dots :\chi_{m_N}(x)),$$
%Denoting    $ N:=\#\big((n-1)\Delta\cap {\mathbb Z}^n\big)$, since $D$ is very ample, we have the embedding $\varphi:X\to{\mathbb P}^N$ determined by the map
%$$\varphi:X\to{\mathbb P}^N,\;\; \varphi(x)=(\chi_{m_1}(x):\dots :\chi_{m_N}(x)),$$
%$$x\in T\mapsto (\chi_{m_1}(x):\dots :\chi_{m_N}(x)),$$
where $m_1,\dots,m_N$ are the points of the set $((n-1)\Delta)\cap {\mathbb Z}^n$.
 Moreover, the inverse image $\varphi^*({\mathcal O}_{{\mathbb P}^N}(1))$ is the line bundle ${\mathcal O}_X(D)$ determined by the divisor $D$. We put $j$ for the composition $Y\hookrightarrow X\overset{\varphi}{\rightarrow}{\mathbb P}^N$. Then ${\mathcal L}$ is the inverse image
\begin{equation}\label{mathcalL}
{\mathcal L}={\mathcal O}_Y\otimes_{j^{-1}({\mathcal O}_{{\mathbb P}^N})}j^{-1}{\mathcal O}_{{\mathbb P}^N}(1).
\end{equation}
% $i$ for the inclusion $i: Y\hookrightarrow X$, and  then $(\varphi\circ i)^*({\mathcal O}_{{\mathbb P}^N}(1))={\mathcal L}$. 

Since $Y$ may be a singular variety, one defines the first Chern class $c_1({\mathcal L})$ 
%can be regarded 
as an  element of the ring Chow $A^*(Y)$ of $Y$; more precisely,    $c_1({\mathcal L})$ is the class in $A^*(Y)$ of  the divisor of any nonzero rational section of ${\mathcal L}$ \cite[page 41]{Fulton-Int} \cite[page 37]{E-H1}.  Thus,   
 if $B$ is a subvariety of $Y$, then  
$c_1({\mathcal L})\cap [B]=c_1({\mathcal L}|_B)$. As ${\mathcal L}={\mathcal O}_{{\mathbb P}^N}(1)|_Y$,
$$c_1({\mathcal L})\cap[B]=[H]\cap[B],$$
  $H$ being a hyperplane divisor in ${\mathbb P}^N$.

 %Through the  ring homomorphism $A^*(X)\to H^*(X,\,{\mathbb Z})$,  from the Chow ring of $X$ to the  cohomology, we will also consider $c_1({\mathcal L})$ as an element of $H^2(Y,\,{\mathbb Z})$.
 
\smallskip

\noindent
 $C^{\infty\,} ${\bf sheaves.}
 %As ${\mathcal L}$ is  ${\mathcal O}_Y$-module and   $\tilde\nabla$ is a connection on ${\mathcal O}_{{\mathbb P}^N}(1)$, which is  an  ${\mathcal O}_{{\mathbb P}^N}$-module, we need to carry out some changes of rings in order to define the connection $\nabla$ on ${\mathcal L}$. 
From the mathematical point of view, the gauge fields are connections on complex vector bundles  over real manifolds. The existence of partitions of the unity in the ${ C}^{\infty}$ category allows to patch local connections \cite{K-N}. So, to study gauge fields on ${\mathcal L}$ we do not need neither the holomorphic structure of ${\mathcal L}$ nor the complex structure of $Y$.

We denote by ${\mathcal A}_{{\mathbb P}^N}$ the sheaf of differentiable functions on ${\mathbb P}^N$. We put ${\mathcal A}_Y$ for
% the pullback ${\mathcal A}_Y:=j^{-1}({\mathcal A}_{{\mathbb P}^N})$; that is,  ${\mathcal A}_Y$ is
 the sheaf on $Y$ which defines the differentiable structure of $Y$. In view of  the above observation, we will consider connections on the $C^{\infty}$ line bundle
\begin{equation}\label{sfL}
      {\sf L}={\mathcal A}_Y\otimes_{j^{-1}({\mathcal A}_{{\mathbb P}^N})}j^{-1}{\mathcal O}_{{\mathbb P}^N}(1).
 \end{equation}
% where   ${\mathcal O}(1):={\mathcal O}_{\mathbb P}(1).$ 

A gauge field on a line bundle over a smooth manifold is given locally by differential $1$-forms. As $Y$ is not necessary smooth, to define gauge fields on $Y$, we will to consider K\"ahler differentials.
                           
 To simplify notations, {\em in this Section}  we
 write ${\mathbb P}:={\mathbb P}^N$,  
%${\mathcal S}:={\mathcal A}_Y$, ${\mathcal T}:={\mathcal A}_{\mathbb P}$ 
 and ${\mathcal R}$ for the pullback $j^{-1}{\mathcal A}_{\mathbb P}$ of ${\mathcal A}_{\mathbb P}$  by the inclusion $j:Y\hookrightarrow{\mathbb P}$. We also put ${\mathcal O}(1):={\mathcal O}_{\mathbb P}(1).$

Given $U$ an open subset of $Y$, $\Omega_{{\mathcal A}_Y(U)/{\mathbb C}}$ denotes the module of 
K\"ahler differentials of ${\mathcal A}_Y(U)$ over ${\mathbb C}$ \cite[Chapter 16]{Eisenbud}.
Let ${\mathcal A}^1_{Y}$ be   the sheaf of ${\mathcal A}_Y$-modules defined by
 $$ {\mathcal A}^1_{Y}(U)   =\Omega_{{\mathcal A}_Y(U)/{\mathbb C}}. $$
 %We also write ${\mathcal A}^1_Y$ for denoting this sheaf of differentials on $Y$.
% obviuosly, $ {\mathcal A}^1_{\mathcal S}={\mathcal A}_Y$.
Similarly, we define the sheaves of K\"ahler differentials
$${\mathcal A}^1_{\mathcal R}(U):=\Omega_{(j^{-1}{\mathcal A}_{\mathbb P})(U)/{\mathbb C}},\;\;\; 
 {\mathcal A}^1_{\mathbb P}(W):=\Omega_{{\mathcal A}_{\mathbb P}(W)/{\mathbb C}}. $$
 %We also put ${\mathcal A}^1_{\mathbb P}$ for denoting the   sheaf ${\mathcal A}^1_{\mathcal T}$.
 
 %The following well-known general fact will be used in the sequel.
 In general,  given a sheaf of rings ${\mathcal C}$ on a space $Z$ and  the ${\mathcal C}$-modules
   ${\mathcal H}$ and ${\mathcal G}$, we denote by $\theta$ the natural morphism  
$$\theta:{\mathcal H}\otimes_{p,{\mathcal C}}{\mathcal G}\to {\mathcal H}\otimes_{\mathcal C}{\mathcal G},$$
between the tensor product presheaf and the tensor product. Since the stalks of ${\mathcal H}\otimes_{p,{\mathcal C}}{\mathcal G}$ and ${\mathcal H}\otimes_{{\mathcal C}}{\mathcal G}$ at any point of $ Z$ are  canonically isomorphic, there exists an open covering of $Z$, such that the restriction of $\theta$ to each element of this covering is   an isomorphism.

 Returning to our case. By the compactness of $Y$, 
 %and according to the recalled fact, 
 $Y$ can be covered by a finite family ${\sf U}=\{U\}$ of ``small enough" open sets
 %$\{U\}$, 
 such that, for any $U\in{\sf U}$
\begin{itemize}
 \item[(a)] The restriction ${\sfç L}|_U$ is a trivial $C^{\infty}$ line bundle.
\item[(b)] The epimorphism  $\pi:{\mathcal R}\to {\mathcal A}_Y$ of sheaves of rings corresponding to the embedding $j$, defines a surjective map 
%$$\pi_U: f\in {\mathcal R}(U)\mapsto f\circ j\in{\mathcal S}(U).$$ 
 $$\pi_U: b\in {\mathcal R}(U)\mapsto \pi_U(b)\in{\mathcal A}_Y(U).$$
\item[(c)] From the following part of exact conormal sequence \cite[page 387]{Eisenbud}, 
$$q:{\mathcal A}_Y\otimes_{\mathcal R}{\mathcal A}^1_{\mathcal R}\to{\mathcal A}^1_{Y}\to 0,$$
one obtains a surjective ${\mathcal A}_Y(U)$-linear map
\begin{equation}\label{(b)}
q_U:{\mathcal A}_Y(U)\otimes_{{\mathcal R}(U)}{\mathcal A}^1_{\mathcal R}(U) \to{\mathcal A}^1_Y(U),\;\;\; c\otimes db\mapsto cd(\pi_U(b)).
\end{equation}
 \end{itemize}
 %We denote by ${\sf E}$ of open subsets $W$ of ${\mathbb P}$ such that $W\cap Y \supset U$.

Let ${\mathcal D}$ be any of the following sheaves on ${\mathbb P}$: ${\mathcal A}_{\mathbb P}$, ${\mathcal O}(1)$, ${\mathcal A}^1_{\mathbb P}$. Given an open subset $U\subset Y$, let ${\sf E}$ denote the  set consisting of  of open  $W\subset{\mathbb P}$ such that $W\cap Y \supset U$.
 Any section of   $j^{-1}{\mathcal D}(U)$ is defined by a family
$$\{ \tilde\zeta_W\in {\mathcal D}(W) \}_{W\in{\sf C}},$$
where ${\sf C}$ is a cofinal subset of  ${\sf E}$, 
% set of open  $W\subset{\mathbb P}$ such that $W\cap Y \supset U$)
 satisfying $\tilde\zeta_W|_{W'}=\tilde\zeta_{W'},$
%$$\tilde\zeta_W|_{W'}=\tilde\zeta_{W'},$$
 for $W'\subset W$. Conversely, such a family defines a section of  $j^{-1}{\mathcal D}(U)$, which will be denoted $j^{-1}(\tilde\zeta)$.

\smallskip

\noindent
{\bf Connections on ${\sf L}$.} 
 A connection on ${\sf L}$ is a morphism $\nabla$ of
 the category $\mathfrak{Sh}$
% ${\mathcal O}_Y$-modules
$$\nabla:{\sf L}\to {\mathcal A}^1_Y\otimes_{{\mathcal A}_Y} {\sf L},$$
satisfying the following property:
% Given point $x\in Y$ and an open $ W\subset Y$ with $x\in W$, there exists an open neighborhood $U$, $x\in U\subset W$, such that for $g\in{\mathcal O}_Y(U)$ and $s\in{\mathcal L}(U)$
Given an open subset $U\subset Y$,  a function $g\in{\mathcal A}_Y(U)$ and $s\in{\sf L}(U)$ then
\begin{equation}\label{nabla_U(g}
\nabla_U(gs)=\theta_U(dg\otimes s)+g\nabla_Us.
 \end{equation}
%where
% the first summand on the right hand side is, in fact the image in the
Here, $\theta_U$ is the canonical map ${\mathcal A}^1_Y(U)\otimes_{{\mathcal A}_Y(U)}{\sf L}(U)\to
({\mathcal A}^1_Y\otimes_{{\mathcal A}_Y} {\sf L})(U).$

If $U$ is a  small enough open  subset  of $Y$, then condition (\ref{nabla_U(g}) implies
\begin{equation}\label{nabla_U(g)b}
\nabla_U(gs)=dg\otimes s+g\nabla_Us.
 \end{equation}

Given a trivialization of ${\sf L}$   on a small enough open $U\subset Y$, a connection on ${\sf L}$ defines a section $\alpha_U\in{\mathcal A}^1_Y(U)$. A general connection on ${\mathcal O}(1)$  gives rise to 
  local $1$-forms on ${\mathbb P}$, but  the restriction to $Y$ of those forms are not necessarily  sections of the sheaf of K\"ahler differentials ${\mathcal A}^1_Y$.

 Nevertheless, some    connections
  % $\tilde\nabla$
   on ${\mathcal O}(1)$ determine  connections on ${\sf L}$.
  % $\nabla$ on ${\sf L}$. 
   Let ${\sf U}$ be the above mentioned open covering of $Y$. 
  %Since ${\sf L}$ is locally free, there is  an open covering ${\sf U}$ of $Y$ consisting of  small enough open subset $U$, such that
    On each $U\in {\sf U}$   there is a local frame $\sigma$ of ${\sf L}$.
 Then $\sigma\in j^{-1}{\mathcal O}(1)(U)$ is defined by a family indexed by a cofinal subset of ${\sf E}$ 
 $$\{\tilde \sigma_W\in{\mathcal O}(1)(W)\}_{W\in{\sf C}},$$
  satisfying the compatibility condition.
  
  Let us assume that $\tilde \nabla$ is a connection on ${\mathcal O}(1)$, such that  {\em there exist sections} $\tilde\alpha_W \in {\mathcal A}^1_{\mathbb P}(W)$, satisfying $\tilde\nabla\tilde\sigma_W=\tilde\alpha_W\otimes\tilde\sigma_W$.  The compatibility condition of the $\tilde\sigma_W$ implies $\tilde\alpha_W|_{W'}=\tilde\alpha_{W'}$, for $W'\subset W$. That is, the $\tilde\alpha_W$ determine an element $j^{-1}(\tilde\alpha)\in{\mathcal A}^1_{\mathcal R}(U)$.
 Thus,
 $$1\otimes j^{-1}(\tilde\alpha)\in {\mathcal A}_Y(U)\otimes_{{\mathcal R}(U)}{\mathcal A}^1_{\mathcal R}(U).$$
 We set
 \begin{equation}\label{alpha_U}
 \alpha_U:=q_U\big( 1\otimes j^{-1}(\tilde\alpha) \big).
 \end{equation}

\begin{Prop}\label{Prop: alphaU=}
The family  $\{\alpha_U\in{\mathcal A}^1_Y(U)\}_{U\in {\sf U}}$  defines a connection on ${\sf L}$.
 \end{Prop}

{\it Proof.} Let $\sigma'$ be other frame of ${\sf L}|_U$. 
%This frame determine the respective sections $\tilde\sigma'_W$ and
Then there exists a function $b\in{\mathcal A}_Y(U)$ different from zero everywhere, such that
$\sigma=b\sigma'$.
The frame $\sigma'$ determines the corresponding sections $\tilde\sigma'_W$ for each
$W$ belonging to a cofinal  subset of ${\sf E}$. Thus,    for each $W$ belonging to a cofinal subset,
% ${\sf C}\subset{\sf E}$,
there is a function ${\tilde b_W}\in{\mathcal A}_{\mathbb P}(W)$ satisfying $\tilde\sigma_W=\tilde b_W\tilde\sigma'_W$ and $\tilde b_W$ is nonzero everywhere. So, one has $j^{-1}(\tilde b)\in {\mathcal R}(U)$ and  
\begin{equation}\label{pUtildeb}
\pi_U(j^{-1}(\tilde b))=b,
\end{equation} 
since $\sigma=\pi_U(j^{-1}(\tilde b))\sigma'$.

 The connection form $\tilde\alpha'_W$ of  $\tilde\nabla$ in the frame $\tilde\sigma'_W$ is defined by the relation $\tilde\nabla\tilde\sigma_W'=\tilde\alpha'_W\otimes\tilde\sigma'_W$. We have the form $\alpha'_U $ defined from the $\tilde\alpha_W'$ as in (\ref{alpha_U})
 \begin{equation}\label{alpha_U'}
 \alpha'_U:=q_U\big( 1\otimes j^{-1}(\tilde\alpha') \big).
 \end{equation}

 Since $\tilde\nabla$ is a connection 
$$\tilde\alpha_W =\tilde\alpha'_W +\tilde b^{-1}_W d\tilde b_W,$$
for $W$ belonging to a cofinal subset of ${\sf E}$. Therefore,
$$j^{-1}(\tilde\alpha)=j^{-1}(\tilde\alpha')+j^{-1}(\tilde b^{-1})j^{-1}(d\tilde b),$$
and
$$q_U(1\otimes j^{-1}(\tilde\alpha))=q_U(1\otimes j^{-1}(\tilde\alpha'))+q_U\big( (1\otimes (j^{-1}(\tilde b))^{-1} d(j^{-1}(\tilde b))) \big).$$ 
From 
%the ${\mathcal A}_Y(U)$-linearity of $q_U$, 
(\ref{alpha_U'}), (\ref{alpha_U})   and (\ref{pUtildeb}), we deduce 
\begin{equation}\label{alphaU=}
\alpha_U=\alpha_U'+b^{-1}db.   
 \end{equation}
Thus, the K\"ahler differentials $\alpha_U$, defined for the open sets of  the covering ${\sf U}$ of $Y$, satisfy the compatibility condition for defining a  connection on ${\sf L}$.
\qed

 \smallskip

\noindent
{\it Remark.} We have constructed  a connection on ${\sf L}$ by restricting a connection on ${\mathcal O}(1)$ whose connection forms, relative to local frames on the members of ${\sf U}$, are K\"ahler differentials.
Conversely,  the connections on ${\sf L}$ can be extended to connections on ${\mathcal O}(1)$.

 Let $\{V_i\},$ $i=0,\dots,N$ be the standard affine covering of ${\mathbb P}$ and let $g_{ij}$ denote the  transition functions of ${\mathcal O}(1)$ corresponding to the usual trivializations on the $V_i$. A connection on ${\sf L}$ is defined by a collection $\alpha_0,\dots, \alpha_N$, with 
$\alpha_k\in{\mathcal A}_Y^1(Y\cap V_k)$, satisfying on $Y\cap V_j\cap V_i$
\begin{equation}\label{alphak=alphaj}
\alpha_i=\alpha_j+g^{-1}_{ji}{ d}g_{ji},
 \end{equation}
since the transition functions of ${\sf L}$ are the restrictions of the corresponding transition functions of ${\mathcal O}(1)$. By (\ref{alphak=alphaj}), together with the fact that $g^{-1}_{ji}{ d}g_{ji}$ is already a $1$-form on $V_j\cap V_k$, there exists for each $j$ an open subset of $W_j\subset V_j$ and 
$\tilde\alpha_j\in{\mathcal A}_{{\mathbb P}^m}^1(W_j)$, such that
$$\tilde\alpha_i=\tilde\alpha_j+g^{-1}_{ji}{ d}g_{ji},$$
on $W_j\cap W_i$. In this way, the forms $\tilde\alpha_j$ define a  connection on ${\mathcal O}(1)|_{V}$, where $V$ is an open subset of ${\mathbb P}$ containing  $Y$. This connection can be extended to a connection  on ${\mathcal O}(1)$ \cite[page 67]{K-N}.

\smallskip

\noindent
{\bf Chern connection on ${\mathcal O}_{{\mathbb P}^N}(1)$.} 
On  ${\mathcal O}(1)$
%${\mathcal O}_{{\mathbb P}^N}(1)$ 
 is defined the Chern connection $\tilde\nabla^{0}$, compatible with the Hermitian and holomorphic structures of this line bundle. 
 The homogeneous coordinates $z_0,\dots,z_N$ of ${\mathbb P}$ define sections of
 %${\mathcal O}_{{\mathbb P}^N}(1)$.
  ${\mathcal O}(1)$.
 We denote by ${\bf v}=(v_1,\dots, v_N)$ the inhomogeneous coordinates on an affine subset $V_k$ belonging to the standard covering of ${\mathbb P}$. The Hermitian metric on 
  ${\mathcal O}(1)|_{V_k}$ is given by the function $h=(1+|{\bf v}|^2)^{-1}$. Thus, the form  of the connection $\tilde\nabla^{0}$ on $V_k$ is
  \begin{equation}\label{equ:tilde_alpha}
  \tilde\alpha_k=h^{-1}\partial h=- (1+|{\bf v}|^2)^{-1} (\bar{{\bf v}}\cdot d{\bf v}).
  \end{equation}
  The curvature of $\tilde\nabla^{0}$ is $\frac{2\pi}{i}\omega_{FS}$, where $\omega_{FS}$ is the Fubini-Study form. Thus,
 \begin{equation}\label{c_1}
 c_1({\mathcal O}_{{\mathbb P}^N}(1))=[w_{FS}]\in H^2({\mathbb P},\,{\mathbb Z}).
  \end{equation}

	 In particular, $\tilde\alpha_k$ {\em is a section of} ${\mathcal A}^1_{\mathbb P}(V_k)$.
	Let ${\sf U}$ be the covering of $Y$ consisting of the open set $U_k=Y\cap V_k$, with $k=0,1,\dots, N$. The form $\tilde\alpha_k$ determines a K\"ahler differential $\alpha_k\in{\mathcal A}^1_Y(U)$, and by Proposition \ref{Prop: alphaU=} the family $\{\alpha_k\}$ defines a connection on ${\sf L}$.

We write $\nabla^0$ for the connection on ${\sf L}$ determined by $\tilde\nabla^{0}$.    If $U$ is an open subset of $Y$ contained in the affine subset $V_k\subset{\mathbb P}$, then 
%of the standard covering of ${\mathbb P}$, 
taking into account (\ref{alpha_U}), (\ref{(b)}) and (\ref{equ:tilde_alpha}),  the connection form of $\nabla^0$ 
can be written on $U\subset Y$ in terms of $\pi_U(v_i)$, the functions on $Y$ induced by the inhomogeneous coordinates $v_i$, as 
\begin{equation}\label{nabla0}
-\pi_U(h)\big((\pi_U\bar{\bf v})\cdot d(\pi_U \bf v)\big).
 \end{equation}
Analogously, for the curvature $d\alpha+\alpha\wedge\alpha$ of $\nabla^0$, one has 
\begin{equation}\label{Fnabla0}
\frac{\pi_U(h^{-1})d\pi_U{\bf v} \wedge d\pi_U\bar{\bf v}-\big(\pi_U\bar{\bf v}\cdot d\pi_U{\bf v}\big)\wedge \big(\pi_U{\bf v}\cdot d\pi_U\bar{\bf v}\big)  }{\pi_U(h^{-2})},
\end{equation}
where the dot product is also involved in $d\pi_U{\bf v} \wedge d\pi_U\bar{\bf v}$. That is, on the smooth locus of $Y$ the curvature $F_{\nabla^0}$ of $\nabla^0$ is 
\begin{equation}\label{Fnabla0}
F_{\nabla^0}=\frac{2\pi}{i}\,\Hat\omega_{FS},
 \end{equation}
 where $\Hat\omega_{FS}$ is the restriction  of the Fubini-Study form to that locus.

 \smallskip
 
\noindent
{\bf Yang-Mills functional.}    We denote by $\Hat g$ the  metric   on  $Y\setminus{\rm Sing}(Y)$  induced by the Fubini-Study one of ${\mathbb P}$.  With $\star$ we denote the corresponding Hodge star operator.
 On the affine space of connections of ${\sf L}$, one defines the Yang-Mills functional
 \begin{equation}\label{mathcalYM}
 \nabla\mapsto{\sf YM}(\nabla)=-\int_YF_{\nabla}\wedge \star F_{\nabla},
  \end{equation}
 where $F_{\nabla}$ is the curvature form of $\nabla$.
 %and   $\star$ is the Hodge star. 

From the variational principle, one deduces that the fields $\nabla$ that are stationary points  for the functional ${\sf YM}$ satisfy the equation
 \begin{equation}\label{YMequation}
\nabla^*F_{\nabla}=0,
\end{equation}
$\nabla^*$ being 
%By the Euler-Lagrange equation, the curvature of the  fields $\nabla$ that are critical points for the Yang-Mills functional satisfy the condition $\nabla^*F_{\nabla}=0$, where $\nabla^* $ is 
the formal adjoint operator of $\nabla$. 
% On the other hand,  Bianchi identity $\nabla F_{\nabla}=0$ holds for any connection.
 As ${\sf L}$ is a line bundle, Bianchi identity and (\ref{YMequation}) reduce to 
%these equations say that
 $$dF_{\nabla}=0,\;\;\;d^*F_{\nabla}=0,$$
where $d^*$ is the codifferential operator.
 That is,  the Yang-Mills fields 
%$\nabla$
 are those whose curvature 
%$F$ 
is a $2$-form closed and coclosed in $Y\setminus{\rm Sing}(Y)$.  

 By (\ref{Fnabla0}),
%On the other hand, the curvature of   $\tilde\nabla^{\rm Chern}$ is $\frac{2\pi}{i}\omega_{FS}$. Thus,
  $${\sf YM}(\nabla^0)=4\pi^2\int_Y  \Hat\omega_{FS}\wedge\star\, \Hat\omega_{FS},$$
% where $\Hat\omega_{FS}$ is the restriction of $\omega_{FS}$ to $Y$. By (\ref{cerradaycocerrada}),  $\Hat\omega_{FS}$ is  Yang-Mills field on $Y$.

 %Denoting by $g$ the Fubini-Study metric and by $\star$ the corresponding Hodge star,
On the other hand,
 \begin{equation}\label{starHodge}
\Hat \omega_{FS}\wedge\star\,\Hat\omega_{FS}=\Hat g(\Hat\omega_{FS},\,\Hat\omega_{FS}){\it vol}={\it vol},
  \end{equation}
 where ${\it vol}$ is the volume form determined by $g$ on the smooth locus of ${Y}$. As $\Hat\omega_{FS}$ and ${\it vol} $  are closed forms,
 $\Hat\omega_{FS}\wedge d(\star\,\Hat\omega_{FG})=0$. Since $\Hat\omega_{FS}$ is non degenerate, it follows $d^*\Hat\omega_{FS}=0$.
%, where $d^*$ is the operator codifferential.
 Thus, $\nabla^0$ is a Yang-Mills field. Furthermore,
 %\begin{equation}\label{cerradaycocerrada}
 %d\Hat\omega_{FS}=0,\;\;\; d^*\Hat\omega_{FS}=0.
 %\end{equation}
 from (\ref{starHodge}), it follows 
 %$$\Hat\omega_{FS}\wedge\star\, \Hat\omega_{FS}=\Hat g(\Hat\omega_{FS},\,\Hat\omega_{FS}){\it vol}={\it vol},$$
% $\Hat g$ being the restriction of the Fubini-Study metric to $Y$ and ${\it vol}$ the volume form on $Y$ defined by $\Hat g$. Thus,
  $$ {\sf YM}(\nabla^0)=4\pi^2{\rm vol}(Y).$$
  
  Similarly, the functional $\sf{YM}_{\mathbb P}$ defined on the connections over ${\mathcal O}(1)$ takes 
  %its minimum the value 
  at the Chern connection $\tilde\nabla^{0}$ the value,
  % and
 $${\sf YM}_{\mathbb P}(\tilde\nabla^{0})=4\pi^2{\rm vol}({\mathbb P}),$$
 and it is the minimum value taken by $ {\sf YM}_{\mathbb P}$.
 
 If $\nabla$ were a gauge field  on $Y$ such that $ {\sf YM }(\nabla)<{\sf YM }(\nabla^0)$,  then one could construct on an small  open neighborhood $V$ of $Y$
% subset $V\supset Y$,
 an extension $\tilde\nabla$ of $\nabla$, such that  
 $${\sf YM }_V(\tilde\nabla^{0})-{\sf YM}_V(\tilde\nabla)=\epsilon>0.$$
 On the other hand,  $\tilde \nabla$ can be extended to ${\mathbb P}$, so that
 $$|{\sf YM }_{{\mathbb P}\setminus V}(\tilde\nabla^{0})-{\sf YM}_{{\mathbb P}\setminus V}(\tilde\nabla)|<\epsilon/2.$$ Thus
  $${\sf YM }_{\mathbb P}(\tilde\nabla)<{\sf YM }_{\mathbb P}(\tilde\nabla^{0}).$$
  But this inequality  contradicts    the fact that ${\sf YM }_{\mathbb P}(\tilde\nabla^{0})$ is the minimum value of ${\sf YM }_{\mathbb P}$.
  Thus, we can summarize the above results in the following   following proposition. 
	%Hence, $\mathcal{YM }(\nabla^{0})$ is the minimum value of the Yang-Mills functional on gauge fields over ${\sf L}$. 
 
 \begin{Prop}\label{Prop:YM} The Yang-Mills functional   
 (\ref{mathcalYM}) reaches its minimum value at the connection $\nabla^0$. 
 \end{Prop}

%%%%%%%%%%%%%%%%%%%%%%%%%%%%%%%%%%%%%%%%%%%%%%%%%%%%%%%%%%%%%%%%%%%%%%%%%%%%%%%%%%%%%%%%%%%%%%%%%%%%%%%%%%%%%%%%%%%%%%%%%%%%%%%%%%%%
%%%%%%%%%%%%%%%%%%%%%%%%%%%%%%%%%%%%%%%%%%%%%%%%%%%%%%%%%%%%%%%%%%%%%%%%%%%%%%%%%%%%%%%%%%%%%%%%%%%%%%%%%%%%%%%%%%%%%%%%%%%%%%%%%%%

\section{Equivariant Chern class}\label{S:Equiv}
 {\bf Universal $T$-bundle.} Let $T$ denote the torus 
 %of the variety $X$; i.e. 
$T=({\mathbb C}^{\times})^n$. The space $({\mathbb P}^{\infty})^n$ can be taken as the classifying  space $BT$ and  the natural 
 $T$-bundle 
$$ET=({\mathbb C}^{\infty}\setminus \{0\})^n\to BT=({\mathbb P}^{\infty})^n,$$
as the  the universal $T$-bundle.
As it is well-known, the rational cohomology $H^*(BT,\,{\mathbb Q})$ is the polynomial ring ${\mathbb Q}[t_1,\dots,t_n]$, where each $t_i$ has degree $2$.

Given $m=(m_1,\dots, m_n)\in {\mathbb Z}^n$, one has the group homomorphism 
 \begin{equation}\label{chi_m}
\chi_m:	g=(z_1,\dots,z_n)\in T\mapsto \prod_i(z_i)^{m_i}\in{\mathbb C}^{\times}.
 \end{equation}
In Appendix, we will prove the following lemma.

\begin{Lem}\label{Lem:Appendix}
The character $\chi_m$ determines: 
\begin{enumerate} 
\item
A map $\xi:BT\to B{\mathbb C}^{\times}$, such that the morphism of algebras induced  by $\xi$ on the cohomologies
 $$\xi^*:{\mathbb Q}[t]=H^*({\mathbb P}^{\infty},\,{\mathbb Q})\to H^*(BT,\,{\mathbb Q})={\mathbb Q}[t_1,\dots,t_n]$$
		satisfies $\xi^*(t)=\sum_{i=1}^n m_i t_i.$								
\item A bundle map $\psi:ET\to E{\mathbb C}^{\times}$ over $\xi$, satisfying
$\psi(e\cdot g)=\chi_m(g)\psi(e)$, for all $e\in ET$ and all $g\in G$. 
 \end{enumerate}
\end{Lem}

		Associated with the $T$-action on the toric variety $X$, there is the homotopy quotient space $ET\times_T X$, whose cohomology is the $T$-equivariant cohomology $H^*_T(X)$ of $X$. The universal $T$-bundle $ET$ does not admit a finite dimensional model. 
		Nevertheless, as  ${\mathbb C}^{\infty}=\cup_r{\mathbb C}^{r}$,    there are finite dimensional obvious approximations   $ET_r\to BT_r$ to $ET \to BT$ that are algebraic varieties. Furthermore,   $H^i(ET_r\times_T X )=H^i_T(X)$	for $r>r(i)$.	Similarly, one defines the $i$-th $T$-equivariant Chow group
		$A_{T,i}(X)$ as the usual Chow group $A_{i+r-n}(ET_r\times_T X )$ \cite{E-G}.
		
		The anticanonical divisor of $X$ is $T$-invariant, hence the invertible sheaf ${\mathcal L}'$ is $T$-equivariant. 
   The $T$-equivariant first Chern class $c_1^T({\mathcal L}')$ is by definition the first Chern class of the line bundle
  $${\mathcal L}'_T:=ET\times_T{\mathcal L}'\to ET\times_TX.$$
		 That is, $c_1({\mathcal L}'_T)$ is the class in $A_{T,1 }(X)$ of divisor defined by any nonzero rational section of ${\mathcal L}'_T$. As $X$ is a rationally smooth variety, then $c_1({\mathcal L}'_T)$ can be considered as a cohomology class in $H^2_T(X,\,{\mathbb Q})$.
		 % by means of the ring homorphism
		% $$A^*(X)\otimes_{\mathbb Z}{\mathbb Q}\to H^*(X,\,{\mathbb Q}),$$
		% where $A^i(X)=A_{n-i}(X)$.

  %That is, $c_1^T({\mathcal L}')=c_1({\mathcal L}'_T)\in H_T^2(X,\,{\mathbb Q})$ \cite{Fulton}. 
  
	\smallskip
	\noindent
	{\bf  $T$-invariant hypersurface.}
  Let us assume that {\it  $Y$ is a $T$-invariant subvariety} of $X$, then ${\mathcal L}$ is a $T$-equivariant line bundle over $Y$, and the inclusion $i:Y\hookrightarrow X$
  induces natural map $\tilde i:ET\times_ T Y\to ET\times_TX$. Over $ET\times_T Y$ we have the line bundles  $\tilde i^*({\mathcal L}'_T)$ and  ${\mathcal L}_T$. 
   As ${\mathcal L}=i^*{\mathcal L'}$, one has   ${\mathcal L}_T=\tilde i^*({\mathcal L}'_T)$. 
  
  By the above equality and the  functoriality of the Chern class 
  \begin{equation}\label{equivChern}
	c_1^T({\mathcal L})=c_1({\mathcal L}_T)=\tilde i^*(c_1^T({\mathcal L'})).
	\end{equation}
  
  The set $X^T$ of fixed point of $X$ for the $T$-action is in bijective correspondence with the set of vertices of the polytope $\Delta$. Since the action of $T$ on $X^T$ is trivial and 
 the rational cohomology of the classifying space 
$ET/T$ is the polynomial ring ${\mathbb Q}[t_1,\dots, t_n]$, the inclusion $X^T\subset X$ gives rise to  the localization map
  $${\lambda}': H^*_T(X,\,{\mathbb Q})\to H^*(ET\times_T X^T,\,{\mathbb Q})=\prod_{v}{\mathbb Q}[t_1,\dots,t_n],$$
  where $v$ ranges the set of vertices of $\Delta$.
  
  The obvious inclusion $Y^T\subset X^T$ gives rise  to the natural projection on the respective equivariant cohomologies
  $$\prod_v{\mathbb Q}[t_1,\dots,t_n]\overset{{\rm pr}}{\longrightarrow}\prod_w{\mathbb Q}[t_1,\dots,t_n],$$
  $w$ ranging on the fixed points in $Y$.
  
Similarly, one has a localization map
$$ \lambda: H^*_T(Y,\,{\mathbb Q})\to H^*(ET\times_T Y^T,\,{\mathbb Q})=\prod_w{\mathbb Q}[t_1,\dots,t_n].$$
 From (\ref{equivChern}), it follows
\begin{equation}\label{eq:pr}
{\rm pr}(\lambda'(c_1^T({\mathcal L}'))=\lambda(c_1^T({\mathcal L})).
 \end{equation}

In what follows, we will determine $c_1^T({\mathcal L}')$ in order to obtain the localization $\lambda(c_1^T({\mathcal L}))$ by means of (\ref{eq:pr}).   $c_1({\mathcal L}')$ is the 
class $[D]\in H^2(X,\,{\mathbb Q})$ of the divisor 
 $\sum_F(n-1)D_F$. This divisor is the one defined for the character
 $$\chi: T\to {\mathbb C}^{\times},\;\;\;\; (z_1,\dots,z_n)\mapsto \prod_i (z_i)^{n-1}.$$
Thus, if $pt$ is any point of the set $D\cap X^T$, the restriction  of ${\mathcal L}'$ to $\{pt\}$ is the vector space $V\simeq {\mathbb C}$ endowed with the $T$-action 
$$g\cdot z=\chi(g)z,$$
with $z\in{\mathbb C}$ and $g\in T$.
So, one can to construct the $T$-equivariant line bundle $V_T=ET\times_TV\to BT$.

 On the other hand, we can endow $V$ with a structure of ${\mathbb C}^{\times}$-equivariant bundle, by means of the natural action of ${\mathbb C}^{\times}$ on
 ${\mathbb C}$. So, one has
 $V_{{\mathbb C}^{\times} } \overset{p}{\to} {\mathbb P}^{\infty},$
  where   $V_{{\mathbb C}^{\times} }$ is the set of classes $[(x_k),\,z]$, with
  $$[(x_k),\,z]= [(\mu x_k),\,\mu^{-1}z],$$
	where $(x_k)\in{\mathbb C}^{\infty}\setminus \{0\}$, $z\in {\mathbb C}$ and  $\mu\in{\mathbb C}^{\times}$.
  Denoting with $\{x_k\}$  the point $[x_0:x_1:\dots]\in{\mathbb P}^{\infty}$, 
 the fiber $p^{-1}(\{x_k\})$ over $\{x_k\}$  is the set
  $$\{[\mu\cdot(x_k),1]\,|\,\mu \in {\mathbb C}^{\times}\}\cup\{0\}.$$
  Thus,  $V_{{\mathbb C}^{\times} }$ is the tautological bundle ${\mathcal O}_{{\mathbb P}^{\infty}}(-1)$.
 %$$V_{{\mathbb C}^{\times} }=\big(({\mathbb C}^{\infty}\setminus \{0})\times {\mathbb C }  \big)/\sim\longrightarrow {\mathbb P}^{\infty}.$$
 The following diagram, where $h$ is defined by $h([e,\, z])=[\psi(e),\,z]$, shows that $V_T$ is the pullback of ${\mathcal O}_{{\mathbb P}^{\infty}}(-1)$ by means of $\xi$.
$$
\xymatrix{
V_T\ar[d]\ar[r]^h & V_{{\mathbb C}^{\times} }\ar[d]^p \\
BT\ar[r]^{\xi} & {\mathbb P}^{\infty} \,.}
$$
From the functoriality of the Chern class together with Lemma \ref{Lem:Appendix}, it follows
$$c_1(V_T)=\xi^*(c_1({\mathcal O}_{{\mathbb P}^{\infty}}(-1)))=-\xi^*(t)=-(n-1)\sum_it_i\in{\mathbb Q}[t_1,\dots, t_n].$$

%We write ${\sf c}$ for denoting the polynomial $-(n-1)\sum_it_i$.
 The preceding arguments are valid for each point of $D\cap X^T$. Thus, the localization $\lambda'$ maps $c_1^T({\mathcal L}')\in H^2_T(X,\,{\mathbb Q})$ to
%$$\big(-\xi^*(t_1,\dots, t_n)\big)^r=\big((1-n)\sum_i t_i\big)^r\in \bigoplus_{{\mathbb P}^N}{\mathbb Q}[t_1,\dots,t_n],$$
$$(\overbrace{{\sf c},\dots,{\sf c}}^{r})\in \prod_v{\mathbb Q}[t_1,\dots,t_n],$$
$r$ being the number of vertices of the polytope $(n-1)\Delta$  and ${\sf c}$ denoting the polynomial $-(n-1)\sum_it_i$.
 %elements of $D\cap X^T$.

From (\ref{eq:pr}), it follows the following theorem. 
\begin{Thm}\label{Th:ewquiv}
Assumed that $Y$ is a $T$-invariant subvariety of $X$, the localization $\lambda(c_1^T({\mathcal L}))$ of the $T$-equivariant Chern class of ${\mathcal L}$ is the element
%$$\lambda(c_1^T({\mathcal L}))= \big((1-n)\sum_i t_i\big)^s\in \bigoplus_w{\mathbb Q}[t_1,\dots,t_n],$$
% $$\lambda(c_1^T({\mathcal L}))= 
$$(\overbrace{{\sf c},\dots,{\sf c}}^{s})\in \prod_w{\mathbb Q}[t_1,\dots,t_n],$$
  %\big((1-n)\sum_i t_i\big)^s\in \bigoplus_w{\mathbb Q}[t_1,\dots,t_n],$$
where $w$ ranges on the set of fixed points of $Y$ and $s$ is the cardinal of this set.
\end{Thm}

%%%%%%%%%%%%%%%%%%%%%%%%%%%%%%%%%%%%%%%%%%%%%%%%%%%%%%%%%%%%%%%%%%%%%%%%%%%%%%%%%%%%%%%%%%%%%%%%%%%%%%%%%%%%%%%%%%%%%%%%%%%%%%%%%%%%%%%%%%%%%%%%%%%%%%%%%%%%%%%%%%%

\section*{Appendix}
In this section we will sketch a proof of Lemma \ref{Lem:Appendix}.
The universal bundle $ET$ is homotopically equivalent to  $EU(1)^n$, the one   of the
 %compact torus with $G$ the 
compact torus $U(1)^n$. According to the Milnor construction 
(see \cite{Mi}, \cite[page 54]{Husemoller}),  given a topological group $G$,  one can take as $EG$ the infinite join 
$$EG=G\ast G\ast G\cdots .$$
 Following Husemoller, we
write $\langle g,t\rangle$ for the elements of $EG$, where
$\langle g,t\rangle$ is given by the
 sequence $( t_0g_0, \, t_1g_1,\, \dots)$ with $t_i\in [0,\,1]$, $g_i\in G$ satisfying the properties detailed  in
 \cite{Husemoller}.
The quotient of $EG$ by the right $G$-action
$$\langle g,\,t\rangle g'=\langle   t_0g_0g', \, t_1g_1g',\, \dots   \rangle$$
is the classifying space $BG$.

Denoting with $EG_r$ the subspace of $EG$ consisting of the elements $\langle g,\,t\rangle$ such that $t_i=0$ for $i>r$, and by $BG_r$ the corresponding quotient, the fibre bundle $EG_r\to BG_r$ is a finite dimensional approximation to the universal bundle.

 When $G=U(1)$,  then $EG_r$ is the sphere $S^{2r+1}$ and $BG_r={\mathbb P}^r$. Moreover, 
$BG_1={\mathbb P}^1 \subset {\mathbb P}^{\infty}=BG$ is the generator of $H_2({\mathbb P}^{\infty},\,{\mathbb Q})$. If $G=U(1)^n$, 
then $BG_1=({\mathbb P}^1)^n \subset ({\mathbb P}^{\infty})^n=BG.$ 

The group homomorphism (\ref{chi_m}) defines a map
$$\psi:=E\chi_m:\langle g,\,t\rangle\in EU(1)^n \mapsto \langle\chi_m(g),\,t\rangle\in EU(1).$$
That is, for $(y_1,\dots,y_n),\,(z_1,\dots,z_n)\in U(1)^n$
$$\psi\langle t_0(y_1,\dots,y_n),\, t_1(z_1,\dots,z_n),\dots\rangle=
\langle t_0\prod_iy_i^{m_i},\,t_1\prod_iz_i^{m_i},\dots\rangle.$$
The induced mapping $\xi:=B\chi_m: BU(1)^n\to BU(1)$ maps
$$\xi\langle t_0(y_1, \,1,\dots,\, 1),\, t_1(z_1,\, 1,\dots,\,1),\,0,\,0,\dots\rangle=
\langle t_0 y_1^{m_1},\,t_1z_1^{m_1},\,0,\,0,\dots\rangle.$$

In other words, $\xi$ restricted to ${\mathbb P}^1\times {\rm pt}\times \dots\times{\rm pt}\subset ({\mathbb P}^{\infty})^n$ is defined by 
%on the first factor of the product ${\mathbb P}^1 \times\dots\times {\mathbb P}^1$ maps
$$\big([y_1:\,z_1],\, {\rm pt},\, \dots ,{\rm pt} \big)\in ({\mathbb P}^1)^n\mapsto [y_1^{m_1} :\, z_1^{m_1}]\in {\mathbb P}^1.$$
%$$\in {\mathbb P}^1\mapsto [z_1^{m_1}:\,y_1^{m_1}]\in {\mathbb P}^1.$$
In terms of the homology,  $\xi$ maps the generator $a_1\in H_2(({\mathbb P}^{\infty})^n,\,{\mathbb Q})$ determined by first factor in the product $({\mathbb P^1})^n$ to $m_1b$, where $b$ is the generator
of $H_2(BU(1),\,{\mathbb Q})$. 

Analogously,
 $$\xi\langle t_0(1,\,\dots,\, 1,\,y_n),\, t_1(1,\,\dots,\, 1,\,z_n),\,0,\,0,\dots\rangle=
\langle t_0 y_n^{m_n},\,t_1z_n^{m_n},\,0,\,0,\dots\rangle.$$
Thus, $\xi$ maps the generator $a_n\in H_2(({\mathbb P}^{\infty})^n,\,{\mathbb Q})$ determined by the $n$-th factor in the product $({\mathbb P^1})^n$ to $m_nb$. 
%, where $b$ is the generator of $H_2(BU(1),\,{\mathbb Q})$. 
That is, 
$$\xi_*:H_*(({\mathbb P}^{\infty})^n,\,{\mathbb Q})={\mathbb Q}[t_1,\dots, t_n]\to 
H_*({\mathbb P}^{\infty},\,{\mathbb Q})={\mathbb Q}[t],$$
maps $t_i$ to $m_i t$.

%In other words, $\xi$ restricted to $({\mathbb P}^1)^n\subset ({\mathbb P}^{\infty})^n$ is defined by 
%on the first factor of the product ${\mathbb P}^1 \times\dots\times {\mathbb P}^1$ maps
%$$\big([z_1:\,y_1],\dots,\,[z_n:\,y_n]\big)\in ({\mathbb P}^1)^n\mapsto [\prod_iz_i^{m_i} :\, \prod_iy_i^{m_i}]\in {\mathbb P}^1.$$
%$$\in {\mathbb P}^1\mapsto [z_1^{m_1}:\,y_1^{m_1}]\in {\mathbb P}^1.$$
%In terms of the homology,  $\xi$ applies the generator $a_i\in H_2(({\mathbb P}^{\infty})^n,\,{\mathbb Q})$ determined by $i$-th factor in the product $({\mathbb P^1})^n$ to $m_ib$, where $b$ is the generator
%of $H_2(BU(1),\,{\mathbb Q})$. 
%Analogously for the other factors in the product $({\mathbb P}^1)^n$. 

The map induced between the cohomologies
$\xi^*:{\mathbb Q}[t]\to{\mathbb Q}[t_1,\,\dots,\,t_n]$
%$$\xi^*:H^*({\mathbb P}^{\infty},\,{\mathbb Q})={\mathbb Q}[t]\to
%H^*(({\mathbb P}^{\infty})^n,\,{\mathbb Q})={\mathbb Q}[t_1,\,\dots,\,t_n],$$
 is the dual of the preceding one, hence
$\xi^*(t)=\sum_i m_it_i$.
\qed

%%%%%%%%%%%%%%%%%%%%%%%%%%%%%%%%%%%%%%%%%%%%%%%%%%%%%%%%%%%%%%%%%%%%%%%%%%%%%%%%%%%%%%%%%%%%%%%%%%%%%%%%%%%%%%%%%%%%%%%%%%%%%%%%%%%%%%
%%%%%%%%%%%%%%%%%%%%%%%%%%%%%%%%%%%%%%%%%%%%%%%%%%%%%%%%%%%%%%%%%%%%%%%%%%%%%%%%%%%%%%%%%%%%%%%%%%%%%%%%%%%%%%%%%%%%%%%%%%%%%%%%%%%%%%%%%

\end{document}